\input amstex
\documentstyle{amsppt}

\magnification=\magstep1
 \def\codim{\operatorname{codim}}
 \def\cal{\Cal}

 \centerline{\bf
  ON THE SEVERI VARIETIES OF SURFACES IN $\bold P^3$}
 \vskip1cm
 \centerline{ L.CHIANTINI }
 \centerline{Universita' di Siena - Dipartimento di Matematica}
 \centerline{Via del Capitano, 15\ \ \ 53100\ SIENA (Italy)}
 \centerline{e-mail: chiantini\@unisi.it}
 \bigskip  
 \centerline{ C.CILIBERTO}
 \centerline{Dipartimento di Matematica -
 Universita' di Roma - Tor Vergata}
 \centerline{ loc. La Romanina, \ \ \ 00133\ ROMA (Italy)}
 \centerline{ e-mail: ciliberto\@axp.mat.utovrm.it}
 
\bigskip
\topmatter
\abstract
For a smooth surface $S$ in $\bold P^3$ of degree $d$ and for 
positive integers $n,\delta$, the Severi variety $V^0_{n,\delta}
(S)$ is the subvariety of the linear system $|O_S(n)|$ which 
parametrizes curves with $\delta$ nodes. We show that for $S$ general, 
$n\ge d$ and for all $\delta$ with $0\le\delta\le\dim(|O_S(n)|)$, 
then  $V^0_{n,\delta}(S)$ has at least one component which is 
reduced, of the expected dimension $\dim(|O_S(n)|)-\delta$. We 
also construct examples of reducible Severi varieties on  general 
surfaces of degree $d\ge 8$.
\endabstract
\endtopmatter

 \bigskip
 
 \heading Introduction \endheading

Nodal curves in the plane have been classically extensively studied. The
variety $V_{n,\delta}$ of irreducible plane curves of degree $n$ with
$\delta$ nodes is called the {\it Severi variety} of type $(n,\delta)$. Severi
claimed in [Se] the irreducibility of $V_{n,\delta}$, but his proof was not
complete. This result has been more recently established by J. Harris in
[H]. \par

It is natural to extend the notion of Severi varieties to curves on
any smooth projective surface $S$. Take an {\it effective} line bundle $\Cal
L$  on $S$, i.e. a line bundle such that the linear system $|\Cal L|$ is not
empty.  Let $p_{\cal L}$ be the arithmetic genus of a curve $C\in |\Cal L|$
and let $\delta$ be an integer such that $0\le \delta\le p_{\cal L}$. It
makes sense to consider the subscheme $V^0_{\Cal L,\delta}$ of the
projective space $|\Cal L|$ consisting of all curves $C\in |\Cal L|$ which
are irreducible and have $\delta$ nodes and no other singular point. This
scheme is locally closed in $|\Cal L|$, and we can consider its closure
$V_{\Cal L,\delta}$ in $|\Cal L|$, which we call the {\it Severi variety} of
curves with $\delta$ nodes in $|\Cal L|$. Often we will also refer to the
open dense subset $V^0_{\Cal L,\delta}$ as to the Severi variety. \par

As for the plane case, one may ask several questions concerning these Severi
varieties. For example:\par

\noindent $(i)$ {\it the existence problem}: when is $V_{\Cal L,\delta}$ not
empty?\par

\noindent $(ii)$ {\it the dimension problem}: what are the dimensions of the
components of $V_{\Cal L,\delta}$?\par

\noindent $(iii)$ {\it the irreducibility problem}: when is $V_{\Cal
L,\delta}$ irreducible?\par

\noindent $(iv)$ {\it the enumeration problem}: what is the degree of
$V_{\Cal L,\delta}$ as a variety in $|\Cal L|$?\par

There is a natural approach to the dimension problem which uses
infinitesimal deformations of nodal curves. Let $C$ be a curve in $V^0_{\Cal
L,\delta}$, which we thus suppose not empty, and let $Z$ be its set of
nodes. One knows (see [S]) that the Zariski tangent space to $V^0_{\Cal
L,\delta}$ at $C$ is \par

$$T_{V^0_{\Cal L,\delta},C}\simeq H^0(S, {\Cal I}_{Z,S}\otimes \Cal L)/(C)$$

\noindent and the corresponding obstruction space is $H^1(S, {\Cal I}_Z\otimes \Cal
L)$.  Let us consider the number $c_{Z,\Cal L}$ of independent conditions
which $Z$ imposes to the linear system $|\Cal L|$, i.e. the rank of the
natural restriction map

$$H^0(S, \Cal L)\to H^0(Z,{\Cal O}_Z)$$

Then one has

$$\dim_C V^0_{\Cal L,\delta}\le \dim  T_{V^0_{\Cal L,\delta},C}=$$

$$=h^0(S, {\Cal I}_Z\otimes \Cal L)-1=dim |\Cal L| - c_{Z,\Cal L}$$

On the other hand one has also 

$$\dim_C V^0_{\Cal L,\delta}\ge h^0(S, {\Cal I}_Z\otimes \Cal L)-1- h^1(S,
{\Cal I}_Z\otimes \Cal L)-1=$$

$$=dim |\Cal L| -\delta$$

\noindent and max$\{dim |\Cal L| -\delta, -1\}$ is the so-called {\it expected
dimension} of $V^0_{\Cal L,\delta}$ at $C$. In conclusion $V^0_{\Cal
L,\delta}$ is smooth of codimension $\delta$ in $|\Cal L|$ at $C$ if and
only if $\delta\le dim |\Cal L| $ and $Z$ imposes independent conditions to
$|\Cal L|$. In this case one says that $V^0_{\Cal L,\delta}$ and $V_{\Cal
L,\delta}$ are {\it regular} at $C$.  If instead $\dim_C V^0_{\Cal
L,\delta}>dim |\Cal L| -\delta$ we will say that $V^0_{\Cal L,\delta}$ is
{\it superabundant} at $C$. \par

Regularity is a very strong property. In fact, again by using infinitesimal
deformation arguments (see [S]), one sees that it implies that {\it the
nodes of $C$ can be independently smoothed}. In other words for any subset
$Z'$ of $Z$, of order $\delta'$ we can find a family of curves in $|\Cal L|$, parametrized by a
disc $\Delta$, whose general member $C_t$, $t\in \Delta-\{0\}$, is nodal with $\delta'$
nodes, its special fibre over $0$ is $C$, and the set of nodes of the general
member specializes to $Z'$ when $C_t$ specializes to $C$.\par

It is easy to see that if $S$ is a rational surface and $C\cdot K_S<0$ for
$C\in |\Cal L| $, then, if $C$ sits in $V^0_{\Cal L,\delta}$, its set of
nodes $Z$ imposes independent conditions to $|\Cal L| $, hence all
components of $V^0_{\Cal L,\delta}$ are regular. This applies in particular
to the plane or, more generally, to the del Pezzo surfaces, for which one
has that all components of the Severi varieties are regular. \par

A similar argument also works for other regular surfaces like K3 or Enriques
surfaces, for which again one proves that all components of the Severi
varieties are regular (see [AC], [T1,2]). \par

Let us now turn to the existence problem. In the plane case one takes
advantage from the fact that there are rational nodal curves of any degree
$n$, namely the general projections to a plane of rational normal curves of
degree $n$ in $\bold P^n$.  By smoothing some of their nodes one proves the
non emptiness of $V^0_{n,\delta}$ for all $\delta \le \frac {(n-1)(n-2)}2$. \par

An investigation about non emptiness of
$V^0_{\Cal L,\delta}$  for del Pezzo surfaces is contained in the recent paper
[GLS]. \par

The existence problem becomes more difficult for other surfaces, like K3 or
Enriques surfaces. For example take an effective line bundle $\Cal L$ on a
K3 surface $S$. Notice that $p_{\Cal L}=dim |\Cal L|$, hence we may expect
the existence of rational nodal curves in $|\Cal L|$. They would be clearly
finitely many and regular, and this in turn would imply that $V^0_{n,\delta}$ is
regular for all integers $\delta$ such that $0\le \delta\le dim |\Cal L|$.\par

Recently Xi Chen proved in [Ch] the following result: {\it for all $g\ge 3$,
on the general K3 surface $S$ of the principal series in $\bold P^g$ there
are rational nodal curves in the linear system $|\Cal O_S(n)|$ for all
$n>0$}. We will comment on this theorem in \S2, where we will recall the
main idea of its proof in the case $g=3$. \par

In recent times there has also been a lot of interest about the enumerative
problem. We recall that the degree of the Severi variety of rational curves
of degree $n$ has been computed by Kontsevich in [K], whereas Caporaso and
Harris gave in [CH] a general formula for the degree of all Severi varieties
in the plane. These formulas were in principle also contained in the paper
[R] by Ziv Ran.\par

We also recall the paper [YZ] by Yau and Zaslow, in which one finds a
suggestion for the computation of the number of rational hyperplane sections
of the general K3 surface $S$ of the principal series in $\bold P^g$. For an
algebro-geometric account of [YZ] one can look at Beauville's paper [B].
Unfortunately there are still some missing details in order to make Yau and
Zaslow's beautiful argument work: among these one should prove that {\it
all} rational hyperplane sections of the general K3 surface $S$ of the
principal series in $\bold P^g$ are nodal. This has been announced, but
still not completely proved, by Chen in [Ch].\par

In the present paper we consider Severi varieties for surfaces of degree $d$
in $\bold P^3$ and in particular we address the existence and dimension
problems. Since the cases $d\le 4$ are covered by the
aforementioned results, we consider the case $d\ge 5$, which is interesting
since the surfaces in question are of general type, and for such
surfaces very little is known about the structure of
their Severi varieties. \par

Let $S$ be a smooth surface of degree $d$ in $\bold P^3$.  We use the
notation $V_{n,\delta}(S)$  and $V^0_{n,\delta}(S)$ to denote the Severi
varieties of curves in $|\Cal O_S(n)|$. A recent result by the first author and
Sernesi [CS] says that {\it all} components of $V^0_{n,\delta}(S)$ are
regular if $\delta$ is {\it small} with respect to $d$ and $n$, for instance
if $\delta<\frac{nd(n-2d+8)}4$. \par

On the other hand, if we are on a surface $S$ of degree $d$, we could in
principle expect the existence of regular components for all $\delta$ such
that $0\le \delta\le dim |\Cal O_S(n)|:=N(d,n)$. This is excactly what we
investigate in the present paper. \par

Our results are of two types. First, in \S 1, we give examples of
superabundant components $V^0_{n,\delta}(S)$ on any surface $S$ of degree
$d$. Our examples are of two kinds. One kind consists of not empty
components of $V^0_{n,\delta}(S)$ with $\delta>N(d,n)$, which are clearly
superabundant. The second kind instead consists of superabundant components
in the range $\delta\le N(d,n)$, $n\ge d$. Both examples are obtained by intersecting cones of a suitable family with $S$. \par

After having recalled in \S 2 Chen's theorem mentioned above, in \S 3 we
extend it to general surfaces of higher degrees. Our main result is theorem
3.1, which says that if $S$ is a general surface of degree $d\ge 5$ in
$\bold P^3$, for any $n>0$ and for any $\delta$ such that $0\le \delta\le
N(d,n)$ there exist regular components of 
$V^0_{n,\delta}(S)$. The proof is by induction, the first step of which is
Chen's theorem. Then we use degeneration and we prove the existence of
surfaces of degree $d$ which split into a general surface of degree $d-1$
plus a general tangent plane of it. We can prove that on such reducible
surfaces we have {\it limit nodal curves} with good properties. By partially
smoothing these curves we find our regular components of the Severi
varieties on a general surface.\par

We finish by noticing that, putting together our main theorem 3.1 and the
examples of \S 1, we conclude with the surprising result that on a general
surface $S$ of large enough degree in $\bold P^3$ some Severi varieties are
reducible. \par

We pose however the following question: {\it is it true that for $S$ general
of degree $d$ and for all $\delta$ such that $0\le \delta\le N(d,n)$, there
is only one regular irreducible component of $V^0_{n,\delta}(S)$?} \par

We also have a question about superabundant components, which naturally
arises by looking at our examples in \S 2: {\it  is it true that for $S$
general all the superabundant components of $V^0_{n,\delta}(S)$ are such
that their general curve is cut out on $S$ by a surface which is singular in
codimension one?}\par

In conclusion, the second author would like to thank J. Harris for several
interesting discussions on the subject of this paper. Both authors wish to thank the referee for useful remarks on the first version of the paper.\par

 \vskip.5cm

\heading 0. Notation \endheading
 
 We work over the complex field $\bold C$. In this paper we will use the
following notation:\bigskip
 
 $S$ will be a general surface of degree $d$ in $\bold P^3$.\par

 $\bold P_d=|\Cal O_{\bold P^3}(d)|$ is the linear system of surfaces of
degree $d$ in  $\bold P^3$. We have:

 $$ \dim\bold P_d=\frac {d^3+6d^2+11d}6.$$\par

 $p_a(d,n)=\frac {dn}2(n+d-4) + 1$ is the arithmetic genus of the divisors
in $|\Cal O_S(n)|$.\par\smallskip

 $V_{n,\delta}(S)$ is the Severi variety of curves in $|\Cal O_S(n)|$ with
$\delta$ nodes, i.e. the closure in $|\Cal O_S(n)|$ of the scheme
$V^0_{n,\delta}(S)$ of irreducible curves with $\delta$ distinct nodes and
no other singularity. Notice that $V^0_{n,\delta}(S)$ is open in
$V_{n,\delta}(S)$.  If there is no danger of confusion, we will write
$V_{n,\delta}$ and $V^0_{n,\delta}$ instead of $V_{n,\delta}(S)$ and
$V^0_{n,\delta}(S)$.\par\smallskip

 $N(d,n)=h^0\Cal O_S(n)-1=\dim|\Cal O_S(n)|$. We have:

 $$N(d,n)=\cases \frac{nd}2(n-d+4)+
 \frac{d^3-6d^2+11d}6-1 & for\ n\ge d \\ \dim\bold P_n & for\ n<d
 \endcases$$
 
We say that $C\subset S$ is a {\it nodal curve} if it is irreducible, with
only nodes as singularities.  $V^0_{n,\delta}$ parametrizes all nodal
curves of geometric genus $p_a(n,d)-\delta$ in the system $|\Cal
O_S(n)|$.\par\medskip
 
Notice that on some dense Zariski open subset $U$ of $\bold P_d$ we have a
universal surface $\Cal S \to U$  and a universal linear system $|\Cal
O_{\Cal S}(n)|$, and one can also define the universal algebraic schemes
$\Cal V_{n,\delta}$ and $\Cal V^0_{n,\delta}$ corresponding to the Severi
varieties.

 \vskip.5cm
 \heading 1. Superabundant Components of Severi Varieties\endheading
 
 As we said in the introduction, in the plane, i.e. when $d=1$, it is
classically known that the Severi varieties $V_{n,\delta}$, when not empty,
all have the expected dimension $N(d,n)-\delta$. This result extends to
every Del Pezzo surface. In particular it applies also to quadrics and
cubics, i.e. to the cases $d=2,3$.\par\smallskip

 On the other hand, when the degree $d$ increases and $S$ is of general
type, the situation is no longer so simple. On a general surface of degree
$d\gg 0 $ the Severi varieties may have several components, with different
dimensions.\par\smallskip

 In particular, the following examples show the existence of superabundant
components of $V^0_{n,\delta}$. Later on we will see that for the values of
$n$, $d$ and $\delta$ considered in example (1.2) below there are also
regular components, hence the corresponding Severi varieties turn out to be
reducible.\par
 
 \proclaim{Example 1.1} \rm Fix $d\ge 20$ and let $S$ be a smooth surface of
degree $d$ in $\bold P^3$. Consider $n=3$, $\delta=20$. We have $N(d,3)=19$
so the expected dimension of $V^0_{n,\delta}$ is negative. \par

 On the other hand, $V^0_{n,\delta}$ is non empty: indeed, take the
intersection of $S$ with a general cone over a singular plane cubic curve.
 \endproclaim
 
 Even when the expected dimension is positive, it may happen that
$V^0_{n,\delta}$ has superabundant components. One can produce such examples
using again cones.

 \proclaim{Example 1.2} \rm Take a smooth surface $S$ of degree $d\ge 8$ and
take $n\gg d$, $n$ even.\par

 On a general plane there is a  family of dimension $(nd-n-2)/2$ of
 irreducible nodal curves of degree $n$ having exactly

 $$\alpha=\frac{n^2-3n+2}2-\frac{nd}2+\frac{7n}2$$

\noindent nodes. Fixing a general $P\in\bold P^3$, we have a family of cones 
with vertex at $P$ over these plane curves. Such cones will have $\alpha$ double
lines. Intersecting these cones with $S$, we get a family contained in
$V^0_{n,\delta}$ with  $\delta=d\alpha$, whose dimension is at least
$(nd-n-2)/2$. In fact two cones cannot intersect $S$ in the same curve $C$,
since $C$ and  $P$ determine the cone.\par

One computes the expected dimension of $V^0_{n,d\alpha}$ which is:

 $$N(d,n)-d\alpha=\frac{d^3-12d^2+11d-6} 6$$

 Notice that it becomes soon smaller than $(nd-n-2)/2$ when $n$ is  large with respect to $d$.
 \endproclaim
 
 Remember from the Introduction that there are no superabundant components of Severi varieties on a quartic surface. In particular, for
$\delta=N(4,n)$, then $V^0_{n,\delta}$ is a finite set of points. There is a
general, easy criterion for the finiteness of $V^0_{n,\delta}$, when $\delta$ is
large.
 
 \proclaim{Proposition 1.4} If $2\delta> n^2d$ then $V^0_{n,\delta}$ is
formed by isolated points.\endproclaim

 \demo{Proof} Assume that $V^0_{n,\delta}$ has a component of positive
dimension. This implies the existence of a non trivial family of nodal
curves, with $\delta$ nodes, in $|\Cal O_S(n)|$. If $C$ is a general element
of this family and $Z$ is its set of nodes, then there exists a non-zero
element in $T_{C,V^0_{n,\delta}}\simeq H^0\Cal I_{Z,S}(n)/(C)$, i.e. there
exists another curve $C'\in|\Cal O_S(n)|$ passing through the nodes of $C$.
Then one has:

 $$2\delta \le C\cdot C' = n^2d.$$\enddemo
 
 \vskip.5cm
 
 \heading 2. Nodal Curves on General Quartic Surfaces \endheading
 
 In this paragraph we consider the case $d=4$. Notice that
$N(4,n)=2n^2+1=p_a(4,n)$ and, if $N(4,n)-\delta\ge 0$, then
every component of $V^0_{n,\delta}$ is regular.\par

 In the case $S$ is general, the existence
problem, i.e.  the non-emptyness of  $V^0_{n,\delta}$, follows from a general
theorem by X. Chen [Ch], whose result, in the present situation, boils down to
the following:\par
 
 \proclaim{Theorem 2.1} {\rm [Ch]}
 Let $S$ be a general quartic surface in $\bold P^3$.
 Then for all $n>0$, the linear system $|\Cal O_S(n)|$
 contains irreducible nodal rational curves. 
 \endproclaim
 
 There are two consequences of Chen's theorem 2.1:\bigskip
 
\noindent $(i)$ $V^0_{n,\delta}\neq\emptyset$ for all $\delta\le
N(4,n)$ (see the Introduction);\par

 \noindent $(ii)$ looking at the universal Severi varieties, the theorem
implies that $\Cal V^0_{n,\delta}$, for $\delta=N(4,n)$, has a component
 which dominates $\bold P_4$, with 0-dimensional general fiber.\bigskip

It is necessary for our purposes to briefly recall the idea of Chen's proof.
\par

Let $S_0=Q_1\cup Q_2$ be a special quartic in $\bold P^3$ formed by the
union of two general quadrics $Q_1$ and $Q_2$.  The quadric $Q_i$, $i=1,2$,
being general, is smooth, therefore it has two rulings of lines, which we
denote by $H^i_j$, $j=1,2$. Moreover the intersection $E$ of $Q_1$ and
$Q_2$ is a smooth elliptic quartic curve in $\bold P^3$. \par

We consider on $S_0$ the curve $C_0$ defined as follows:\par

$$C_0=C_0^1\cup C_0^2$$

\noindent with $C_0^i$ sitting on $Q_i$, $i=1,2$, defined as:\par

$$C_0^i = C^i_1\cup ... \cup C^i_n$$

\noindent where $C^i_j\in |H^i_1|$, $j=1,...,n-1$, are distinct lines and
$C^i_n\in |H^i_1+nH^i_2|$ is a smooth irreducible curve, subject to the
following conditions:\par

\noindent $(i)$ the line $C^1_j$, $j=1,...,n-1$, intersects $E$ at two
distinct points $q_{2(j-1)}, q_{2j-1}$;\par

\noindent $(ii)$ the line $C^2_j$, $j=1,...,n-1$, intersects $E$ at the two
distinct points $q_{2j-1}, q_{2j}$;\par

\noindent $(iii)$ the curve $C^i_n$ passes through $q_{2(n-1)}$ and has a
contact of order $2n+1$ with $E$ at point $r$ different from all the points
$q_h$, $h=0,...,2(n-1)$;\par

\noindent $(iv)$ the curve $C^2_n$ passes through $q_0$ and has a contact of
order $2n+1$ with $E$ at $r$.\par

Notice that the linear systems $|H^i_1+nH^i_2|$, $i=1,2$, cut on $E$ a
complete linear series of degree $2n+2$ and dimension $2n+1$. Taking this
into account, after one moment of reflection the reader will see that such a
configuration exists, and that $C_0$ is a divisor in $|\Cal O_{S_0}(n)|$.  \par

We will refer to $C_0$ as to a {\it Chen's curve}.\par

Now one takes a general quartic surface $S$ passing through $q_0$, and
therefore cutting $E$ in fifteen more distinct points $p_1, ..., p_{15}$
different from the points $q_h$, $h=1,...,2(n-1)$, and $r$. One takes the
pencil generated by $S_0$ and $S$, which we briefly denote by $S_t=tS+S_0$.
General facts about deformations of singular curves (see [Ch], [R]), ensure
that, for $t$ in a small disc $\Delta$ around $0$, $C_0$ deforms to an
irreducible nodal curve $C_t$ on $S_t$ in the linear system 
$|\Cal O_{S_t}(n)|$, having $2n^2+1$ distinct nodes, and therefore being
rational.  The limits of the nodes of $C_t$, when $t$ tends to $0$ are:\par

\noindent $(i)$ the $n(n-1)$ nodes of $C_0^1$ arising from the intersections
of $C^1_1,...,C^1_{n-1}$ with $C^1_n$;\par

\noindent $(ii)$ the $n(n-1)$ nodes of $C_0^2$ arising from the
intersections of $C^2_1,...,C^2_{n-1}$ with $C^2_n$;\par

\noindent $(iii)$ the point $q_0$;\par

\noindent $(iv)$ the $2n+1$-tacnode $r$ where $2n$ nodes of $C_t$ coalesce.\bigskip

Let now $S$ be a smooth surface in $\bold P^3$ and let us take a curve $C$ of
degree $m$ on $S$. Let $\pi$ be the general tangent plane to $S$. Notice that
$\pi$ is not tangent to $C$, since the dual variety $V^*$ of a given projective
variety $V$ determines $V$. Then the divisor $\Gamma$ cut out by $\pi$ on
$C$ is formed by $m$ distinct points. Let $U$ be the not empty open subset
$U$ of $S^*$ such that for all $\pi\in U$ this happens. We can consider the
incidence correspondence\par

$$X=\{(x,\pi)\in C\times U: x\in \pi\}$$

\noindent The projection $p: X\to U$ is a covering of degree $m$, whose monodromy
group, which is a subgroup of the full symmetric group $S_{m}$ on $m$
objects, we denote by $G(C,S)$.\par

The following result will be an important ingredient for our construction of
regular components of Severi varieties with the maximum number $N(d,n)$ of
nodes.\par 

\proclaim{Theorem 2.2}  If $S$ is a general quartic surface, for all $n$ and
$\delta \leq N(4,n)$ there is an irreducible component of
$V^0_{n,\delta}(S)$ such that for $C$ general in this component one has
$G(C,S)=S_{4n}$.
 \endproclaim

\demo{Proof} Since any irreducible nodal rational curve on $S$ is the limit
of nodal curves with $\delta\leq N(4,n)$ nodes, it clearly suffices to prove
that on a general quartic surface $S$ there is some rational nodal curve $C$
in $|\Cal O_S(n)|$ for which $G(C,S)=S_{4n}$. \par

We consider the setting of Chen's proof of theorem (2.1) and we keep the
notation introduced above. We will prove that a rational curve $C_t$ on
$S_t$, $t\in \Delta-\{0\}$, which tends to a Chen's curve, is such that
$G:=G(C_t, S_t)$, i.e. the monodromy group of the corresponding cover $p_t:
X_t\to U_t$, is the full symmetric group $S_{4n}$.  In order to do so, we
will prove that $G$ contains all transpositions.\par

 First we prove that $G$
contains a simple transposition. To prove this, it suffices to prove that the
monodromy group $G_0=G(C_0,S_0)$ of the limit covering of $p_t: X_t\to U_t$
contains a simple transposition. When $S_t$ tends to $S_0$, the dual variety
$S_t^*$ breaks up in several components, among which there are the dual planes
corresponding to the points $p_i$, $i=1,...,15$. Actually we may assume that
$p_1$ is a general point of $E$. Then the general tangent plane $\sigma$ to
$C^2_n$ through $p_1$ is simply tangent to $C^2_n$ at some point $w$ and
meets transversally the curve $C_0$ at the remaining intersection points.
The existence of such a simply tangent plane yields a transposition in $G_0$
as proved in [ACGH], pg. 112.  Indeed a plane $\pi_0$ passing through $p_1$
and sufficiently close to $\sigma$ cuts $C_0$ transversally, and two points
of the intersection $x_0,y_0$ come together at $w$ when $\pi_0$ tends to
$\sigma$. These two points are exchanged by a simple trasposition in $G_0$.\par

Let us prove now that $G$ contains all transpositions.  Let us take $C_t$
close to $C_0$ and $\pi_t$ tangent to $S_t$ close to $\pi_0$. Let us take
$a_t,b_t$ two points in $\Gamma_t={\pi}_t\cap C_t$. Suppose that the limits
of these two points, when $t$ tends to $0$, are two points $a_0,b_0$ of
$C^2_n$. Then we easily find a permutation in $G_0$ which sends the pair
$(a_0,b_0)$ to the pair $(x_0,y_0)$, simply by moving $\pi_0$ inside the
dual plane to $p_1$, and by letting $a_0$ tend to $x_0$ and $b_0$ to $y_0$.
Since we can transpose the pair $(x_0,y_0)$ with $G_0$, we can also
transpose the pair $(a_0,b_0)$ with $G_0$. This clearly implies that we can
transpose the pair $(a_t,b_t)$ with $G$. \par

In order to conclude, it suffices to prove that, in the same situation as above,
given any pair $(c_t,d_t)$ of points in $\Gamma_t$, there is a permutation in
$G$ which sends the pair $(c_t,d_t)$ to a pair $(a_t,b_t)$ such that the points
of its limiting pair $(a_0,b_0)$ both sit on $C^2_n$. We prove that
acting with $G$ we can move $c_t$ to a point $a_t$ such that its limit $a_0$
sits in $C^2_n$. The same proof, applied to the stabilizer of $a_t$, will show
that we can also move $d_t$ to $b_t$ as above. \par

Suppose $c_t$ tends to a point $c_0$ on a line $C^i_j$, say on $C^1_1$. The
lines $C^1_1$ and $C^2_1$ meet at the point $q_1$ and, locally around $q_1$, we
have an open, smooth neighborhood $A_t$ of $C_t$ which tends to the union
of two neighborhoods of $q_1$ on $C^1_1$ and $C^2_1$ respectively, thus acquiring
a node at $q_1$. By moving $\pi_t$ on $S^*_t$ close to $\pi_0$, we can move the
point $c_t$ describing the whole open subset $A_t$ of $C_t$. Therefore, by
moving $\pi_t$ in this way, we can let the limit point $c_0$ of $c_t$ land on
the line $C^2_1$ instead of $C^1_1$. Our proof is thus finished. \enddemo

 \proclaim{Problem 2.3} \rm  Theorem 2.2 suggests the following
general problem, which would be interesting to investigate. Let $S$ be a smooth
surface of degree $d$ in $\bold P^3$ and let $C$ be an irreducible curve of
degree $m$ on it. Under which conditions is the group $G(C,S)$ equal to $S_m$?
Notice that $G(C,S)$ is certainly not equal to $S_m$ if $S$ is a smooth quadric.
In this case the section of any curve $C$ on $S$ with a tangent plane
to $S$ is formed by two sets of collinear points and the monodromy does not
exchange points of the two sets. A similar situation takes place on any ruled
surface $S$ of degree $d$, but notice that then $S$ is singular as soon as
$d\geq 3$.
\endproclaim
 
\vskip.5cm
 
 \heading 3. Existence of Regular Components \endheading

Our objective in this section is to prove our main result:\par

\proclaim{Theorem 3.1}  Let $S$ be a general surface in $\bold P_d$, $d\geq
4$. For all $n\ge d$ and $\delta \leq N(d,n)$ there is an irreducible
regular component of $V_{n,\delta}(S)$. \endproclaim

As we know from \S 1, it suffices to prove the theorem for $\delta=N(d,n)$,
which means that we have to exhibit an irreducible curve $C$ in $|\Cal
O_S(n)|$ with $N(d,n)$ nodes and no other singularities, and with the nodes
imposing independent conditions to 
$|\Cal O_S(n)|$. Such a curve $C$ gives rise to an isolated, reduced point
of $V^0_{n,N(d,n)}$. We will say that such a curve $C$ is {\it
isolated}.  We will start from Chen's theorem and prove the existence of
these isolated curves inductively. Actually, we will inductively prove the
following statement:\par

\proclaim{Theorem 3.2}  Let $S$ be a general surface in $\bold P_d$, $d\geq
4$. For all $n\ge d$ and there are isolated curves $C$ in  $V^0_{n,N(d,n)}$
such that $G(C,S)$ is the full symmetric group $S_{nd}$. \endproclaim

Our strategy for the proof of theorem 3.2 is the following. We consider a
degenerate surface $S_0$ in $\bold P_{d}$ given as $S_0=F\cup \pi$ where
$F\in \bold P_{d-1}$ and $\pi$ is a plane. We inductively require the
existence of a surface $G\in \bold P_n$ such that $C_0^1=F\cap G$ is
irreducible, nodal and isolated with $N(d-1,n)$ nodes. Furthermore we want
$C_0^2=\pi\cap G$ nodal, with $(n-d+2)(n-d+3)/2$ nodes. Moreover we need
that $C_0^1$ and $C_0^2$ intersect transversally at $n(d-1)$ points of
$F\cap \pi$. Then we consider the curve $C_0=C_0^1\cup C_0^2=S_0\cap G$ and
we see that $C_0$ deforms on a general surface $S$ of degree $d$ to an
irreducible, $N(d,n)$-nodal curve $C$. Using induction we can prove
that $C$ is isolated and that $G(C,S)$ is the full symmetric group. However the
main point is to start by constucting a {\it suitably good} surface $S_0$. Let us
see what the construction of $S_0$ amounts to.\par

First of all we make a remark. If theorem 3.2 holds, then there is a not
empty open subset $U$ of $\bold P_d$ over which the isolated curves of the
statement fit together, giving at least one irreducible component $V_d$ of
the universal Severi variety $\Cal V^0_{n,N(d,n)}$, which dominates $\bold
P_d$ with 0-dimensional general fibers.\medskip

We need some more  notation. \par

Fix numbers $n,d$ with $n+1\ge d\ge 5$.  Assume that theorem 3.2 holds for
$d-1$ and fix a component $V_{d-1}$ of the universal Severi variety
$V^0_{n,N(d-1,n)}$, dominating $\bold P_{d-1}$ with finite fiber over the
general point of $\bold P_{d-1}$.\par
 
Identifying a point $G\in\bold P_n$ with the surface in $\bold P^3$
parametrized by $G$,  we let:\bigskip

 $\bar X\subset\bold P_n\times\bold P_{d-1}$ be a component of the variety
of pairs $(G,F)$ such that $G\cap F$ is an irreducible nodal curve in the
fiber over $F$ of the component $V_{d-1}$ of the universal Severi
variety;\smallskip

 $X$ be the projection of $\bar X$ to $\bold P_n$;\smallskip

 $\bar Y\subset \bold P_n\times(\bold P^3)^*$ be the variety of pairs
$(G,\pi)$ with $G\cap\pi$ irreducible nodal plane curve with
 $(n-d+2)(n-d+3)/2$ nodes;\smallskip

 $Y$ be projection of $\bar Y$ to $\bold P_n$;\smallskip

 $W\subset \bold P_n$ be the variety which parametrizes all reducible
surfaces of degree $n$ which split into a union $F\cup Q$ where $F$ has
degree $d-1$, $Q$ has degree $n-d+1$ and both are irreducible;\smallskip

 \def\PP{\tilde\bold P} $\PP$ be the blow-up of $\bold P_n$ along the
closure of $W$;\smallskip

 $E$ be the exceptional divisor of this blow-up;\smallskip

 $\tilde X$, $\tilde Y$ be the strict transforms of $X$, $Y$; \smallskip

 $\tilde W$ be $\tilde X\cap E$.\bigskip

As indicated in the rough sketch of proof of theorem 3.2 we gave above, the
surface $S_0$ we need corresponds, with our notation, to a point of $X\cap
Y$. In order to prove the existence of {\it suitably good} such surfaces, we
will use degeneration to further reducible surfaces, and we will actually
work on $\PP$ . \par

First we need a few preliminary results. \par
 
 \proclaim{Proposition 3.3} (a) $W$ is irreducible, of dimension

 $$\dim W=\dim \bold P_{d-1} + \dim \bold P_{n-d+1}=$$

$$=2+\frac 16(n^3-3n^2d+3nd^2+9n^2+12d^2-18nd+26n-24d).$$\par

 (b) The tangent space to $W$ at a general point $F\cup Q$ is
 naturally isomorphic to $H^0\Cal I_\gamma(n)/(F\cup Q)$
 where $\gamma $ is the curve $F\cap Q$.\par

 (c) If theorem 3.2 holds for $d-1$ then:\par

$$\dim \bar X=\dim W+1.$$\endproclaim

 \demo{Proof} (a) is trivial. (b) see \S 2 of [CHM]: a particular case is
treated there, but the argument runs in general in the same way.\par

In order to compute the dimension of $\bar X$, project it to $\bold P_{d-1}$
and look at the fibers. On $F\in \bold P_{d-1}$ general, there are only
finitely many nodal curves in $|\Cal O_F(n)|$ with $N(d-1,n)$ nodes,
so the fiber over $F$ is a finite union of projective spaces $\bold P_{n-d+1}$.
Since $\bar X$ dominates $\bold P_{d-1}$, the assertion clearly follows.\enddemo
 
 \proclaim{Proposition 3.4} $Y$ is irreducible, of codimension

 $$\codim_{\bold P_{n}} (Y)=\frac {(n-d)(n-d+5)}2.$$

 \endproclaim

 \demo{Proof} By [H], on a general plane $\pi$, the Severi variety
$V:=V(n,(n-d+2)(n-d+3)/2)$ is irreducible of codimension $(n-d+2)(n-d+3)/2$
in $|\Cal O_\pi(n)|$. $\bar Y$ dominates $V$, with general fibers which are
open subsets of projective spaces of fixed dimension, hence $\bar Y$ is
irreducible and its dimension can be easily computed.\par

By the assumptions, we have $(n-d+2)(n-d+3)\ge 3$ and the general element of $V$
is a nodal plane curve whose nodes are not aligned. Then the projection $\bar
Y\to Y$ has finite fibers. In fact given a smooth surface $G$ of degree
$n\ge 2$ and given a nodal plane section $C$ of $G$ with at least three
nodes which are not aligned, then clearly $C$ is an isolated curve of the
corresponding Severi variety. \enddemo
 
 \proclaim{Proposition 3.5} Assume theorem 3.2 holds for $d-1$. Then $W$
is a proper subvariety of codimension $1$ in the closure of $X$.\endproclaim

 \demo{Proof} Consider a general element $F\cup Q\in W$. Since theorem 3.2
holds for $d-1$, there is some $G\in X$ such that $(F,G)\in \bar X$. Taking
$F$ general, we may assume that also $G$ is general in $X$.\par

 Consider a general element $G'$ of the linear system spanned by $F\cup
Q$ and $G$. Since $G\cap F=G'\cap F$ and $G$ is general in $X$, then
$(F,G')$ belongs to $\bar X$, so that $G'\in X$. Thus in the limit, we get
$F\cup Q\in X$. This proves that $W$ is a proper subvariety of the closure of
$X$.\par 

Now, by proposition 3.3, (c), we have:

$$\dim W< \dim X\leq \dim \bar X=\dim W+1$$

\noindent which concludes the proof of the proposition. 
\enddemo
 
 \proclaim{Remark 3.6} \rm With a similar argument, we may show that a
general element of $X$ is smooth. Indeed, if $(F,G)\in \bar X$ is general,
just as observed above, a general element in the linear system generated by
$G$ and by surfaces of type $F\cup Q$ must belong to $X$.  The assertion
easily follows by applying Bertini's theorem.\endproclaim
 
 \proclaim{Proposition 3.7} For $F\cup Q\in W$ general, the fiber of the
exceptional divisor $E$ over $F\cup Q$ is canonically isomorphic
 to $\bold P(H^0\Cal O_\gamma(n))$, where $\gamma$ is the curve $F\cap Q$.
\endproclaim

 \demo{Proof} The proof runs as in proposition 1 of [CHM] where a particular
case is treated. \enddemo
 
 \proclaim{Proposition 3.8} Suppose theorem 3.2 holds for $d-1$. Then the
general fiber of the projection $\tilde W\to W$ is finite.\endproclaim

 \demo{Proof} Since $W$ belongs to the closure of $X$, the projection
 $\tilde W\to W$ is dominant. By proposition 3.2

 $$\dim W =\dim X-1=\dim \tilde W  $$

\noindent because $X$ is not contained in the exceptional divisor. \enddemo
 
 \proclaim{Remark 3.9} \rm We can say more about the  fibers
 of the projection $\tilde W\to W$.\par

 Take general $F\in\bold P_{d-1}$, $Q\in \bold P_{n-d+1}$ and take $G\in
\bold P_n$ such that $(F,G)$ is a general point of $\bar X$. The pencil
spanned by $F\cup Q$ and $G$ determines, as we saw in proposition 3.5, a
tangent vector to $X$ which is also a normal vector to $W$ in $\bold P_n$,
by proposition 3.3 (b). It follows that this pencil determines an element
$\psi$ of the fiber of $\tilde W$ over $F\cup Q$.\par

 Now consider $\gamma=F\cap Q$ and think of $\psi$ as an element of the
fiber of $E$ over $F\cup Q$, which is identified with
 $\bold P(H^0\Cal O_\gamma(n))$. By construction,  
 $\psi$ corresponds to the divisor $\gamma\cap G$ on $\gamma$. In other words,
for $F$ general, divisors in $|\Cal O_\gamma(n)|$, corresponding to elements of
$\tilde W$ are obtained by intersecting $\gamma$ with the curves in $|\Cal
O_F(n)|$ having $N(d-1,n)$ nodes.\endproclaim 
 
 \proclaim{Lemma 3.10} {\rm [Hi]} Let $\Cal C$ be a flat family of curves,
whose general element is rational and irreducible. Let $C_0$ be a special,
possibly reducible, member of the family. Then every irreducible component
of $C_0$ is rational.\endproclaim
 
 \proclaim{Lemma 3.11}  Let $F$ be general in $\bold P_{d-1}$, $d>4$, and
let $C\in |\Cal O_F(n)|$ be an isolated nodal curve with $N(d-1,n)$ nodes,
such that $G(C,S)$ is the full symmetric group. \par 

Take a general plane $\pi$ tangent to $F$ and take a set $Z$ of $3(n-d+1)-1$
points of $\Gamma=C\cap\pi$. Then there exists a finite, positive number of
rational nodal curves of degree $n-d+1$ on $\pi$, with nodes outside $F$,
passing through $Z$.  More precisely, for each one of these rational nodal
curves the union of the set of nodes and of the set $Z$ gives independent
conditions to the curves of degree $n-d+1$.\endproclaim

 \demo{Proof} The family of rational nodal curves of degree $n-d+1$ in $\pi$
is irreducible of dimension $3(n-d+1)-1$. Thus, if we fix a set $Z$ of
$3(n-d+1)-1$ points in $\Gamma$,  we do get an element, in the closure of this family, which passes through $Z$. We need such an element to be nodal, irreducible, isolated, with its nodes off $F$.\par

Let $\Cal C$ be any family of curves on $\pi$ which is {\it rationally
determined}, i.e.  it is defined over the field of rational functions of
$\pi$. The condition on $G(C,S)$ implies that $\Gamma$ is in uniform position
with respect to any rationally determined family $\Cal C$ of curves on $\pi$, hence
$Z$ and
$\Gamma$ impose the same number of conditions to $\Cal C$, unless $Z$ imposes
independent conditions to $\Cal C$. By Hironaka's lemma 3.2, a degeneration of
rational curves cannot contain the curve $F\cap\pi$, since we may assume that
$F\cap\pi$ has positive genus. It follows that $Z$ must impose independent
conditions to the family $\Cal C$ of  rational curves in  $|\Cal
O_\pi(n-d+1)|$ and to every subfamily of $\Cal C$ which is rationally
determined on $\pi$.\par

This implies that:\par

\noindent $(i)$ $Z$ is contained in a finite number of rational curves
$D\in |\Cal O_\pi(n-d+1)|$;\par

 \noindent $(ii)$ $D$ cannot have singularities worse than nodes, since
curves with such singularities fill a proper, closed subfamily of $\Cal C$; \par

\noindent $(iii)$ similarly, we get that the set $Z'$ of nodes of $D$ does
not intersect $F\cap\pi$;\par

\noindent $(iv)$ $Z\cup Z'$ impose independent conditions to curves of
degree $n-d+1$.\enddemo

We are now in a position to prove the existence of our surface $S_0=F\cup \pi$.
 
 \proclaim{Lemma 3.12} Assume theorem 3.2 holds for $d-1$. Then there is
an irreducible component $B$ of $\tilde Y\cap\tilde X$ not contained in $E$,
whose general point $G$, which we may identify  with its image in $\bold
P{_n}$, has the following properties:\par

 $(i)$ $G$ is irreducible;\par

 $(ii)$ there is a plane $\pi$ such that $\pi\cap G$ is irreducible, nodal,
with $(n-d+2)(n-d+3)/2$ nodes which impose independent conditions to the
curves of degree $n-d+1$, i.e. there are no such curves through these nodes;\par

 $(iii)$ there is a surface $F\in\bold P_{d-1}$ with $F\cap G$ irreducible,
nodal and rigid, with $N(d-1,n)$ nodes;\par

 $(iv)$ the nodes of $F\cap G$ and $\pi\cap G$ are distinct and $G$ is
transversal to $F\cap\pi$.\par

 Furthermore, by moving $G$ in $B$, we may assume $F$ general in $\bold
P_{d-1}$.\endproclaim

 \demo{Proof} Look at $\tilde W\cap \tilde Y$. This intersection has a
component $A$ whose general element is $F\cup Q$, obtained by taking
generically:\par

 - a surface $F\in\bold P_{d-1}$;\par

 - a  plane $\pi$ tangent to $F$;\par

 - a nodal curve $C=F\cap G'$ for some $G'$ with $(F,G')\in\bar X$;\par

 - a subset $Z\subset \Gamma= C\cap\pi$ of degree $3(n-d+1)-1$;\par

 - a nodal rational curve $D\in|\Cal O_\pi(n-d+1)|$, containing $Z$, with
nodes off $F\cap\pi$ (see lemma 3.11);\par

 - a surface $Q$ of degree $n-d+1$, such that $Q\cap \pi=D$.\par

 In fact in this way we found reducible surfaces $F\cup Q$ with a
$N(d-1,n)$-tangent surface $G'$ of degree $n$, which defines the element of
$\tilde W$ associated to $\psi=G'\cap C$ (see remark 3.8). Notice that $\pi$
is tangent to $F$, tangent to $Q$ at $(n-d)(n-d-1)/2$ distinct points and
passes through a subset $Z$ formed by $3(n-d+1)-1$ points of $ \psi$. Then,
by the results of [Ch] and [R], 
the curve $\pi\cap(F\cup Q)$ is the limit of a nodal plane curve $D$ with

 $$ 1+\frac {(n-d)(n-d-1)}2+3(n-d+1)-1 = \frac{(n-d+2)(n-d+3)}2$$

\noindent nodes lying on some irreducible surface $G$ of degree $n$ which
degenerates to $F\cup Q$.\par

 Let us compute the dimension of $A$. We have $\dim \bold P_{d-1}$
parameters for the choice of $F$, two parameters for the choice of $\pi$, no
parameters for the choice of $C$, $Z$ and $D$ and 
 $\dim \bold P_{n-d}$  parameters for the choice of the surface $Q$ passing
through the plane curve $D$. Thus:

 $$\multline
 \dim A=2+\frac 16(n^3-3n^2d+3nd^2+6n^2-9nd+9d^2+11n-9d)=\\
 = \dim W-\codim Y= \dim \tilde X+\dim\tilde Y-\dim\bold P_n -1.
 \endmultline$$

 It follows that there is a component $B$ of $\tilde X\cap\tilde Y$ which
contains $A$ and is not contained in $E$. The general point $G$ of $B$ must
be irreducible, since the elements of $A$ split in two irreducible
components and $G$ cannot split in the same way. The plane $\pi$ and the
surface $F$ exist and fulfill the transversality conditions $(iv)$, since
this happens for general points in $A$. \par

 Finally observe that, by lemma 3.10, there are no curves of degree $n-d+1$
on $\pi$, other than $\pi\cap Q$, passing through $Z$ and through the set
nodes of $\pi\cap Q$. Since $\pi\cap Q$ misses the point where $\pi$ is
tangent to $F$, the nodes of the curve $\pi\cap(F\cup Q)$ also impose
independent conditions to the curves of degree $n-d+1$, i.e. there are no
such curves thorugh these nodes. \enddemo
 
 Now we have all the tools for accomplishing the proof of theorem 3.2.
 
 \demo{Proof of theorem 3.2} The theorem holds for $d=4$, as we saw in \S 2.
We argue by induction and we assume the theorem holds for $d-1$. \par

Take $G$ as in the previous theorem, and take accordingly the surface $F\in
\bold P_{d-1}$ and the plane $\pi$ as  described there. The surface
$S_0=F\cup\pi$ is a reducible surface of degree $d$ and it contains two
curves $C_0^1=F\cap G$ and $C_0^2=\pi\cap G$ which are both nodal, with a
total of

 $$N(d-1,n)+(n-d+2)(n-d+3)/2=N(d,n)$$

\noindent nodes sitting off $F\cap\pi$. These two curves match together at
$n(d-1)$ points of the intersection $F\cap\pi$, where their intersection is
transverse. \par

$C_0^1$ is isolated on $F$ by construction, since we may assume $F$ general
in $\bold P_{d-1}$. $C_0^2$ cannot move on $\pi$ mantaining the number of
nodes and the matching with $C_0^1$: this would imply the existence on $\pi$
of a 1-parameter family of curves of degree $n$, passing through $C_0^2\cap
F$ and through the nodes of $C_0^2$ and this, in turn, would imply the
existence of a curve of degree $n-d+1$ passing through the nodes of
$C_0^2$, which is excluded by lemma 3.12. This implies that the curve
$C_0=C_0^1\cup C_0^2$ is isolated.\par

Now we take a general surface $S$ of degree $d$ and deform $S_0$ in the
pencil generated by $S_0$ and $S$, which, as usual we denote by
$S_t=tS+S_0$. By the results of [Ch] and [R] we see that, for $t$ in a small
disc around $0$, there is on the surface $S_t$ a nodal curve $C_t$,  with
$N(d,n)$ nodes, degenerating to $C_0$ as $t$ goes to $0$. Since $C_0$ is
isolated, the same is true for $C_t$.\par

As for the final assertion of theorem 3.2, its proof can easily be obtained
by imitating the proof of theorem 2.2 and it is even simpler than that.
Therefore we leave it to the reader.\enddemo 
 
 \proclaim{Example 3.13} \rm  Combining theorem 3.1 with example 1.2, one
has examples of reducible Severi varieties $V^0_{n,\delta}$ for general
surfaces of sufficiently high degree in $\bold P^3$.\endproclaim

 \proclaim{Remark 3.14} \rm   If one wants to avoid the complication of
proving theorem 2.2 together with the corresponding inductive assertion on
$G(C,S)$ in theorem 3.2, one can simply take in lemma 3.11 the plane $\pi$ to be
a general plane and not a general tangent plane to $F$. Since the monodromy group
of the general hyperplane section of an irreducible curve is the full symmetric
group (see [ACGH]), one can proceed as in the proof of lemma 3.11. In this
way however one loses a node at a time in the inductive process and therefore
one can only prove the existence of regular components of $V_{n,\delta}(S)$, for
$S$ a general surfaces in $\bold P_d$, $d\ge 4$, for all $n\ge d$ and $\delta \le
N(d,n)-d+4$. In any event this suffices to prove the existence of reducible
Severi varieties on $S$, as indicated in the previous example 3.13.\endproclaim

\vskip.5cm
 
 \Refs

[AC] E. Arbarello, M. Cornalba,  Footnotes to a paper of Beniamino Segre,
Math. Ann. 256 (1981), 341-362\par

[ACGH] E. Arbarello, M. Cornalba, Ph. Griffiths, J. Harris, Geometry of
algebraic curves, Vol. I, Springer Verlag, Berlin, 1985\par

[B] A. Beauville, Counting rational curves on $K3$ surfaces, pre-print, 1997
(Alg. Geom. $97 01 019$) \par

[CH] L. Caporaso, J. Harris, Degrees of Severi varieties, pre-print,
1996 (Alg. Geom. $96 08 025$)\par
 
[Chen] Xi Chen, Rational curves on K3 surfaces, Thesis, Harvad University,
1997\par

[CS] L. Chiantini, E. Sernesi, Nodal curves on surfaces of general type,
Math. Ann., 307 (1997), 41-56\par

[CHM] C. Ciliberto, J. Harris, R. Miranda, General components of the
Noether-Lefschetz locus and their density in the space of all surfaces,
Math. Ann., 282 (1988), 667-680\par

[GLS] G.M. Greuel, Ch. Lossen, E. Shustin, Geometry of families of nodal curves
on the blown-up projective plane, pre-pint, 1997 (Alg. Geom. $97 04
010$)\par
 
 [H] J. Harris, On the Severi problem, Inventiones Math., 84 (1986), 445-461\par
 
 [Hi] H. Hironaka, On the arithmetic genera and the effective genera of
algebraic curves, Mem. College Sci. Univ. Kyoto (A) 30 (1957), 177-195\par

[K] M. Kontsevich, Enumeration of rational curves via torus action,
pre-print\par
 
 [R] Z. Ran, Enumerative geometry of singular plane curves, Inventiones
Math. 97 (1989), 447-465 \par
 
[S] E. Sernesi, On the existence of certain families of curves, Inventiones
Math. 75 (1984), 25-57\par
 
 [Se] F. Severi, Vorlesungen ueber algebraische Geometrie, Teubner, Leipzig,
1921\par

[T1] A. Tannenbaum, Families of algebraic curves with nodes, Compositio Math.
41 (1980), 107-126\par

[T2] A. Tannenbaum, Families of curves with nodes on $K3$ surfaces, Math.
Ann. 260 (1982), 239-253\par

[YZ] S-T Yau, E. Zaslow, BPS states, string duality, and nodal curves on K3,
pre-print, 1996\par
 
 \endRefs
 
 \end